\documentclass[11pt,twoside]{article}
\usepackage{graphics,latexsym,amssymb,amsmath}
\RequirePackage{graphics,epsfig,psfrag}%,picins}
\usepackage{hyperref}
\hypersetup{
    bookmarks=true,         % show bookmarks bar?
    unicode=false,          % non-Latin characters in Acrobat's bookmarks
    pdftoolbar=true,        % show Acrobat's toolbar?
    pdfmenubar=false,        % show Acrobat's menu?
    pdffitwindow=false,     % window fit to page when opened
    pdfstartview={FitH},    % fits the width of the page to the window
    pdftitle={Every knot is a billiard knot},    % title
    pdfauthor={PVK, DP},     % author
    pdfsubject={MEGA11},   % subject of the document
    pdfcreator={PVK},   % creator of the document
    pdfproducer={Producer}, % producer of the document
    pdfkeywords={keyword1} {key2} {key3}, % list of keywords
    pdfnewwindow=true,      % links in new window
    colorlinks=true,       % false: boxed links; true: colored links
    linkcolor=black,          % color of internal links
    citecolor=black,        % color of links to bibliography
    filecolor=black,      % color of file links
    urlcolor=black           % color of external links
}
\def\QQ{{\bf Q}} %reelle Zahlen
\def\pn{\medskip\par\noindent}
\def\Frac#1#2{{\displaystyle{{#1}\over{#2}}}}

\def\[#1\]{\begin{eqnarray}#1\end{eqnarray}}
\def\$#1\${\begin{eqnarray}#1\end{eqnarray}}

\def\phi{\varphi}

\def\eps{\varepsilon}

\def\Lim{\mathop{\hbox{lim}}\limits}
\newcommand{\Pf}{{\em Proof}. }

\newcommand{\EPf}{\hbox{}\hfill$\Box$\vspace{.5cm}}

\def\Frac#1#2{{\displaystyle{\frac{#1}{#2}}}}
\def\phi{\varphi}

\newtheorem{theorem}{Theorem}
\newtheorem{definition}[theorem]{Definition}
\newtheorem{remark}[theorem]{Remark}
\newtheorem{lemma}[theorem]{Lemma}
\newtheorem{proposition}[theorem]{Proposition}

\newtheorem{thm*}[theorem]{Theorem}
\date{\today}
\begin{document}
\pagestyle{myheadings}
\markboth{P. -V. Koseleff, D. Pecker, \today}{{\em Every knot is a billiard knot}}
%%%%%%%%%%%%%%%%%%%%% Publisher's Area please ignore %%%%%%%%%%%%%%
%wis%\catchline{}{}{}{}{}
%%%%%%%%%%%%%%%%%%%%%%%%%%%%%%%%%%%%%%%%%%%%%%%%%%%%%%%%%%%%%%%%%%%
\title{Every knot is a billiard knot}
\author{P. -V. Koseleff \& D. Pecker}
\maketitle
\begin{abstract}
We show that every knot can be realized as a billiard  trajectory  in a convex prism. This solves a conjecture of Jones and Przytycki.
\pn {\bf keywords:} {Billiard knots ,
Lissajous knots, Chebyshev knots, Cylinder  knots}
\pn
\end{abstract}
\section{Introduction}
The study of billiard trajectories in a polyhedron was introduced in 1913 by
K\"{o}nig and Sz\"{u}cs  in \cite{KS}.
They proved density results for a billiard trajectory in a cube.
Their theorem is  strongly related to the famous Kronecker density theorem (see \cite{HW}).

More recently, Jones and Przytycki
  considered the periodic billiard trajectories with no self-intersection as knots.
They proved that billiard knots in a cube are isotopic to Lissajous knots, and
deduced that not all knots are billiard knots in a cube
(\cite{JP}, see also \cite{La,C,BHJS,BDHZ}).
They also proved that every torus knot
(or link) of type $(n,k)$, where $n \ge 2k+1$ can be realized as a billiard
 knot in a cylinder (or in a prism with a regular $n$-gonal floor).
Przytycki went deeper into the study of symmetrical billiards in \cite{P}.

Lamm and Obermeyer \cite {LO} proved that billiard knots in a cylinder
are either periodic or ribbon, hence not all knots are billiard knots in a cylinder.
In \cite{KP1} we constructed many other examples of billiard knots in
convex polyhedrons (in fact irregular truncated cubes).
Dehornoy constructed in \cite{D} a billiard which contains all knots, but this billiard is not convex.
\pn
In this paper we prove the following conjecture of Jones and Przytycki:
\pn
{\em Every knot is a billiard knot in some convex polyhedron.}
\pn
Our result is more precise:\\
{\bf Theorem \ref{th:billiardkl}}
{\em Every knot (or link) is a billiard knot (or link)  in some convex right prism.}
\pn
Using a theorem of Manturov \cite{M}, we first prove that every knot has a diagram which is a star polygon.
Then, perturbing this polygon, we obtain an irregular diagram of the same knot.
We deduce that it is possible to suppose that
$1$ and the arc lengths
of the crossing points are linearly independent over $ \QQ.$
Then, it is possible to use the classical Kronecker density theorem to prove our result.

\section{Every knot has a projection which is a star polygon}\label{star}
A toric braid is a braid  corresponding to the closed braid obtained by projecting
the standardly embedded torus knot into the $xy$-plane.
 A toric braid is a braid of the form
 $ \tau_{k , n} =  \Bigl(\sigma _1 \, \sigma _2 \, \cdots \sigma_{k-1} \Bigr)^n$,
where $ \sigma _1, \ldots, \sigma _{k-1}$ are the standard generators
of the full braid group $B_k.$
A quasitoric braid of type $ (k,n)$
 is a braid obtained by changing some crossings in the toric braid
 $\tau_{k , n}.$
The quasitoric braids form a subgroup of $ B_k.$ Consequently there exist trivial
quasitoric braids of arbitrarily great length, and any quasitoric braid is equivalent to a quasitoric braid of type $(k,n)$  with $ n \ge 2k+1.$

Manturov's theorem tells us that every knot (or link) is realized as the closure of
a quasitoric braid (\cite{M}).

The following definition of polygonal stars will be useful for links.

\begin{definition}Let $p,q$ be integers. The  polygonal star $\{\frac pq \} \subset {\bf R}^2 $
is given by its
vertices $ e(k) = e ^{ \frac {2ki \pi }p},$  and its sides
$ (e(k), e(k+q)),$     $\  k=0, \ldots, p-1  .$
\end{definition}

\medskip

When $p$ and $q$ are coprime integers, this is the usual definition of  star polygons.
The following picture
shows the  polygonal stars $ \{ \frac {10}3 \} ,  \   \{  \frac {10}2  \}$
 and $  \{ \frac 93 \} ,$ as projections of
billiard torus links  $ T(10,3), \ T (10,2)  $ and $T(9,3).$ The dotted lines correspond
to the parts $ z<0$ of the link.

\begin{figure}[th]
\begin{center}
\begin{tabular}{ccc}
\epsfig{file=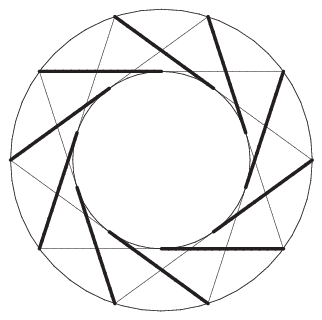,height=3.5cm,width=3.5cm}&
\epsfig{file=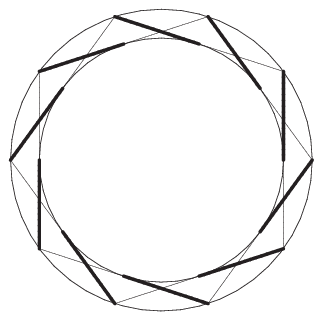,height=3.5cm,width=3.5cm}&%$T(10,3)$&$T(9,3)$
\epsfig{file=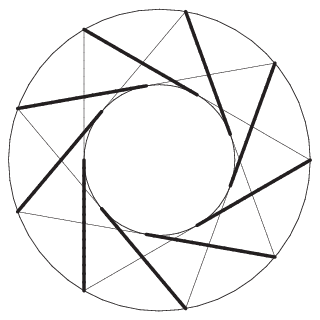,height=3.5cm,width=3.5cm}
\end{tabular}
\end{center}
\caption{The polygonal stars $ \{ \frac {10}3 \}$, $ \{ \frac {10}2 \}$  and $ \{ \frac {9}3 \}$, projections of the torus links
$T(10,3)$, $T(10,2)$ and $T(9,3)$.  }  \label{sp}
\end{figure}

\begin{theorem}\label{th:star}
Every  knot (or link) has a projection that is a  polygonal star.
\end{theorem}
\Pf
Let our  knot be realized as the closure of a quasitoric braid of type $(k,n).$
By our remark, we can suppose $n\ge 2k+1.$
Now, we use the result of Jones and Przytycki which says that every torus knot of
type $ (k,n) , \  n \ge 2k+1 $
can be realized as a billiard knot in a cylinder.
In their construction, the projection on the $xy$-plane is the closure of the toric braid
$ \tau _{ k , n},$ which proves our result.
\EPf

\begin{remark} It is also possible to use a theorem of Lamm and Obermeyer \cite{LO}
 to give another proof of  theorem 1 in the knot case.   A Rosette braid is a braid of the form
 $ \Bigl(\sigma _1 ^{\eps _1} \cdots \sigma _{k-1} ^{ \eps _ {k-1} } \Bigr) ^n
 , \  \eps _i = \pm 1,$
and a Rosette knot is the closure of a Rosette braid.
The theorem of Lamm and Obermeyer tells us
  that every Rosette knot can be represented by a billiard knot in a cylinder.
The knot diagrams obtained in their proof are  star polygons
isotopic to  the  closures of some quasitoric braids.
\end{remark}

\section{Breaking the symmetry}\label{breaking}

Since we want to obtain all knots, we need  irregular diagrams.

First, let us recall some facts about   billiard trajectories.
If $ABC$ is a piece of a polygonal line, then  the mirror placed at  $B$,
is the hyperplane $\mu(B)$ which is orthogonal to the internal bisector of $\widehat{B}$ at  $B$.
The mirror room at $B$ is the closed half-space containing $A, B ,C$ and the mirror at $B.$

We define a billiard trajectory to be a finite union
of polygonal lines, which is  contained in all its mirror rooms.
A  billiard knot (or link) is a polygonal knot (or link)  (\cite{A,C})
which is a  billiard trajectory.

The following result allows us to forget about the billiard, and focus our attention
 on the trajectory. It is valid in every dimension.
\begin{lemma}
Let $ {\cal Q}= (Q_0, \ldots , Q_{n-1}) ,  $  be a billiard trajectory
such that $ {\cal Q} \bigcap \mu (Q_k) = Q_k.$
Then, if $ {\cal P}= (P_0, \ldots , P_{n-1})$ is sufficiently close to
$ {\cal Q},$ it is a billiard trajectory in some convex polyhedron.
\end{lemma}
\Pf
Let $ \overrightarrow{ u _k}({\cal Q}) $ be the unit vector of the internal bisector
of $ \widehat{Q_k}.$
The hypothesis means that for every $k,i, \  i\ne k ,$ the scalar product
$ \overrightarrow{u_k}({\cal Q}) \, {\bf .} \,  \overrightarrow{Q_k Q_i} $ is positive.

Since $\overrightarrow{u_k}({\cal Q}) $  and  $ \overrightarrow{Q_k Q_i}$ depend
continuously on $ {\cal Q},$ this condition remains true for any
trajectory ${\cal P}$ that is sufficiently close to $ {\cal Q}.$
\EPf

\begin{proposition}\label{prop:billiardknot}
Let $K$ be a knot. There exists a plane billiard trajectory  ${\cal P}$
which is a projection
of a knot isotopic to $K,$  and which satisfies the following irregularity condition:

If $ t_i$ are the arc lengths corresponding to the crossings of ${\cal P}$, then
the numbers $ 1, \, t_i$ are linearly independent over $ \QQ .$
\end{proposition}
\Pf
Let $ { \cal Q} = (Q_0, \ldots , Q_{n-1}), Q_n=Q_0$
be a star polygon which is a projection of $K.$
Let us suppose that each line $ (Q_k Q_{k+1})$ has an equation of
the form $ y= \alpha _k x + \beta _k . $
Then, if $(a_k,b_k)$ are sufficiently close to $(\alpha _k, \beta _k),$
the lines
$\{ y= a_kx + b_k \}$ determine a nonconvex polygon
$ {\cal P }= (P_0, \ldots, P_{n-1})  $ close to $ {\cal Q }.$
By our lemma, $ {\cal P}$ is a periodic billiard trajectory in some convex polygon.

By Baire's theorem, we can suppose that the numbers $ a_0,a_1, \ldots, a_{n-1}$ and
$ b_0, b_1, \ldots , b_{n-1}$ are {\bf algebraically} independent over $ \bf Q$.

Let $I$ be the set of integer pairs $ (i,j), \ j \ne i-1,$ such that
the intersection of
 $ [P_i, P_{i+1}]$ and $[P_j, P_{j+1}] $ is a point $ P_{ i,j}$.
\psfrag{P0}{\Large $P_0$}
\psfrag{P1}{\Large $P_1$}
\psfrag{P2}{\Large $P_2$}
\psfrag{P3}{\Large $P_3$}
\psfrag{P4}{\Large $P_4$}
\psfrag{P02}{\Large $P_{0,2}$}
\psfrag{P03}{\Large $P_{0,3}$}
\psfrag{P13}{\Large $P_{1,3}$}
\psfrag{P14}{\Large $P_{1,4}$}
\psfrag{P24}{\Large $P_{2,4}$}
%\psfrag{Pi}{\Large $P_i$}
%\psfrag{Pj}{\Large $P_{j+1}$}
%\psfrag{Pip}{\Large $P_{i+1}$}
%\psfrag{Pjp}{\Large $P_{j}$}
%\psfrag{Pim}{}
%\psfrag{Pipp}{}
%\psfrag{Pjpp}{}
%\psfrag{Pjm}{}
%\psfrag{Pij}{}
%\psfrag{a}{\Large $\ell_{i,j}$}
%\psfrag{b}{}
%\psfrag{c}{}
%\psfrag{d}{}
%\psfrag{e}{}
%\psfrag{f}{}
%\psfrag{b}{$a_2$}%
%\psfrag{c}{\small $a_{n-1}$}\psfrag{d}{\small $a_{n}$}%
\begin{figure}[th]
\begin{center}
{\scalebox{.6}{\includegraphics{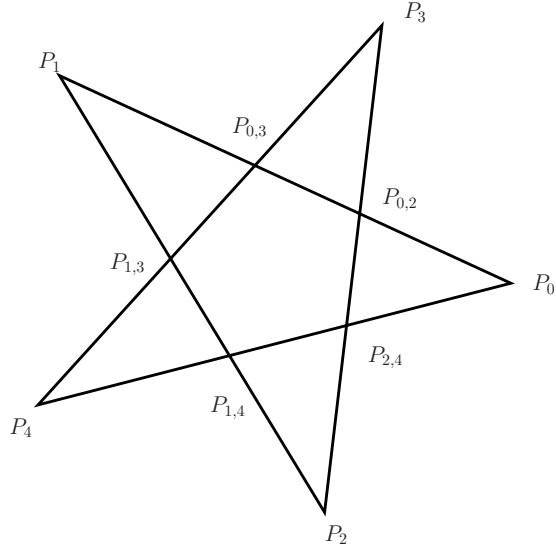}}}%\quad {\scalebox{.6}{\includegraphics{sp2b.eps}}}
\end{center}
\caption{Naming the vertices and crossings of a pentagonal trajectory.}  \label{sp2}
\end{figure}

The vertex $P_i$ is the point $ P_{i-1 , i} .$
The abcissa of $P_{i,j} $  is $ x_{i,j } = \Frac { b_i-b_j}{a_j-a_i},$
and the  length of $[ P_i, P_{i,j} ] $ is $ |\ell_{i,j}|,$  where
$
 \ell_{i,j} =
 \sqrt {1+ a_i^2 }
\Bigl( x_{i,j}   -x_{i-1,i}  \Bigr).
$

 Let us show that the numbers $\ell_{i,j} $ and $1$ are linearly independent
over $ \QQ .$
Suppose that we have a linear relation with rational coefficients
$ \sum_{ (i,j) \in I } \lambda _{i,j} \ell _{i,j} = \lambda$,  with
 $   \lambda , \  \lambda_{i,j}  \in \QQ.$

 This is an algebraic relation between the $a_i$ and the $b_i.$
Since these numbers are algebraically independent over $ \QQ,$ this relation must
be an identity.

Let $k\le n-1$  be a fixed non-negative integer, and let us substitute
$ b_k=1$ and $b_i=0$ if $i \ne k$ in this identity.

%the terms $x_k, \ x_{k+1} , $ and $ x_{k,j}$ are the only ones depending on $b_k.  \  $
 We obtain a new identity between the $ a_i.$
\[
\sum _{ j=0 } ^{n-1} \lambda_{k,j} \sqrt{1+ a_k^2}
\bigl(\Frac 1{a_j-a_k} - \Frac 1 {a_{k-1}-a_k} \bigr)
- \bigl(\sum_ { j=0 }^{n-1}  \lambda_{ k+1,j} \bigr)
\Frac { \sqrt {1+ a_{k+1}^2}} {a_{k+1} - a_k}
&&\nonumber\\
+ \sum _{ i=0 }^{n-1}  \lambda_{i,k}  \Frac {\sqrt {1+a_i^2} }{ a_i-a_k}
&=& \lambda \label{form}
\]
where
$\lambda _{k,k}= \lambda _{ k, k-1}= \lambda _{k+1,k} = 0.$
Substituting $a_k= \sqrt {-1}$, in this identity, we obtain
$$
- \bigl(\sum_ { j=0 }^{n-1}  \lambda_{ k+1,j} \bigr)
\Frac { \sqrt {1+ a_{k+1}^2}} {a_{k+1} - \sqrt {-1}}
+ \sum _{ i=0 }^{n-1}  \lambda_{i,k}  \Frac {\sqrt {1+a_i^2} }{ a_i-\sqrt{-1}} = \lambda.
$$
Let $h\ne k+1$ be an integer,  and let $a_{h} \rightarrow \sqrt {-1}.$
From $\Lim_{z\rightarrow \sqrt {-1}} \Frac {\sqrt { 1+z^2} } {z- \sqrt {-1} }= \infty$, we obtain $\lambda_ {h,k} =0.$
 Since this is true for every $h$ and $k,$ we deduce that
 $ \lambda_{i,j}=0 $ for all $(i,j)\in I$.

Finally, since the arc lengths of the points $P_{i,j} $ are given by
$ t_{i,j}= |\ell_{0,1}| + |\ell_{1,2}| + \cdots + |\ell_{i-1, i}| + |\ell_{i,j}|, $
we deduce the result.
\EPf

\begin{proposition}
Let $L$ be a link. There exists a plane billiard trajectory  ${\cal P}$
which is a projection
of a link isotopic to $L,$  and which satisfies the following  condition.

\smallskip

If $ {\cal R}$ is a component of ${\cal P}$ parametrized by arc length, and
if $ t_i$ are the arc lengths corresponding to the crossings, then
the numbers $ 1, \, t_i$ are linearly independent over $ \QQ .$
\end{proposition}
\Pf
The proof is almost identical to the preceding one.
\pn
There is a link isotopic to $L$ whose plane projection
is a union of polygons
$$ {\cal P } =
 {\cal P}^{ (1)} \cup  {\cal P}^{ (2)} \cup \cdots \cup  {\cal P}^{(d)}$$
whose vertices are  $(P_0, P_1, \ldots , P_{N-1})$.
Let ${\cal R}$ be a component of ${\cal P},$   we can suppose that
the vertices of  ${\cal R}$ are $(P_0, P_1, \ldots , P_{n-1})$.

Furthermore, we can suppose that the equations
$ y = a_k x + b_k $ of the sides of ${\cal P}$ are such that the numbers
$ a_k$ and $b_k, \  k= 1 , \ldots  , N-1, $ are algebraically independent over $ \QQ.$

Here, we consider   the set $I$  of integer pairs
$(i,j), \    i \in \{ 0, \ldots , n-1 \},  \ j \in \{ 0, \ldots , N-1 \}, \  j \ne i-1,$
corresponding bijectively to the arc lengths of the vertices and crossings contained
in $ {\cal R}.$

Then, the rest of the proof is exactly the same as in the case of knots.
\EPf

\section{Proof of the theorem }

We will  use Kronecker's theorem (\cite[Theorem~443]{HW}): %, p. 382):
\begin{theorem}[Kronecker  (1884)]
If $ \theta _1, \theta_2, \ldots, \theta_k, 1$ are linearly independent over $\QQ,$
then the set of points
$ \Bigl((n \theta _1), \ldots , (n \theta_k) \Bigr) $ is dense in the unit cube.
Here $ (x)$ denotes the fractional part of $x.$
\end{theorem}
Now, we can prove our main theorem.
\begin{theorem}\label{th:billiardkl}
Every knot (or link) is a billiard knot (or link) in some convex prism.
\end{theorem}
\Pf
First, we consider knots.
By Theorem \ref{th:star} there exists a knot  isotopic to $K$
whose projection on the $xy$-plane is a periodic billiard trajectory
in a convex polygon $ {\bf D}.$
If $t_i$ are the arc lengths corresponding to the crossings, we can
suppose by Proposition \ref{prop:billiardknot} that the numbers $ t_1, \ldots , t_k, 1$ are linearly
independent over $ \QQ$.
Using a dilatation, we can suppose that the total length of the trajectory is $1.$

Consider the  polygonal curve  defined by $ (x(t), y(t),z(t)),$
where $z(t)$ is the sawtooth function $z(t)= 2 | (nt+ \phi)-1/2 |$
 depending on the integer $n$ and on the real number $\phi.$
If the heights $z(P_k)$ of the vertices are such that $z(P_k) \ne 0, 1$, then
it is a periodic billiard trajectory in the prism $  {\bf D} \times [0,1] $
(see \cite{JP,La,LO,P,KP1}).
If we set $ \phi = 1/2 + z_0/2,$ $ z_0 \in ]0,1[,$ we have $ z(0)=z_0.$
Now, using Kronecker's theorem, there exists an integer $n$ such that
the numbers $z(t_i)$ are arbitrarily close to any chosen collection of heights,
which implies the result.
\pn
The case of links is similar. First, we  find a plane billiard diagram of our link,
and then we  parameterize each component.
\EPf

\begin{remark}
If the diagram has some regularity, then it is generally impossible to use
Kronecker's theorem.
This is illustrated by Lissajous knots and cylinder knots.
This is also true for more general diagrams.

For example, suppose that there are four crossings of arc lengths $ t_1, t_2, t_3, t_4$
such that
$ t_2-t_1= t_4-t_3.$
So, if $z(t)= 2 \bigl | (nt + \phi) -\frac 12 \bigr | $ is the height function, then
$ \eps _i z_i = 2( nt_i + \phi)-1, $ with $   \eps_i = \pm 1.$
We deduce that
$ \eps_1 z_1 - \eps _2 z_2 - \eps _ 3 z_3  + \eps _4 z_4 \in 2 {\bf Z}  .$
Consequently,  we see that
$z_1=z_2=z_3=1$ implies $z_4=1.$
This clearly shows that the heights of the crossings cannot be chosen arbitrarily.
\end{remark}

\vfill
\hrule
\pn
P. -V. Koseleff, UPMC-Paris 6, IMJ and INRIA Paris-Rocquencourt, {\tt koseleff@math.jussieu.fr}\\
D. Pecker, UPMC-Paris 6, Mathematics, {\tt pecker@math.jussieu.fr}
\end{document}